\documentclass[11pt]{article}

\usepackage{psfig,amssymb,latexsym}

\newtheorem{theorem}{Theorem}[section]

\newtheorem{lemma}{Lemma}[section]

\def\slfrac#1#2{\hbox{\kern.1em %
 \raise.5ex\hbox{\the\scriptfont0 #1}\kern-.11em %
 /\kern-.15em\lower.25ex\hbox{\the\scriptfont0 #2}}}

\newcommand{\eqn}[1]{(\ref{#1})}

\newcommand{\eeq}{\end{equation}}
\newcommand{\beql}[1]{\begin{equation}\label{#1}}
\newcommand{\bsq}{{\vrule height .9ex width .8ex depth -.1ex }}

\newcommand{\agh}{{a}}

\newcommand{\ZZ}{{\mathbb Z}}

\newcommand{\QQ}{{\mathbb Q}}
\newcommand{\CC}{{\mathbb C}}

\newcommand{\bU}{{\bf U}}

\newcommand{\sD}{{\cal D}}

\newcommand{\sH}{{\cal H}}

\newcommand{\sR}{{\cal R}}

\makeatletter
\def\@sect#1#2#3#4#5#6[#7]#8{\ifnum #2>\c@secnumdepth
     \def\@svsec{}\else
     \refstepcounter{#1}\edef\@svsec{\csname the#1\endcsname.\hskip .75em }\fi
     \@tempskipa #5\relax
      \ifdim \@tempskipa>\z@
        \begingroup #6\relax
          \@hangfrom{\hskip #3\relax\@svsec}{\interlinepenalty \@M #8\par}%
        \endgroup
       \csname #1mark\endcsname{#7}\addcontentsline
         {toc}{#1}{\ifnum #2>\c@secnumdepth \else
                      \protect\numberline{\csname the#1\endcsname}\fi
                    #7}\else
        \def\@svsechd{#6\hskip #3\@svsec #8\csname #1mark\endcsname
                      {#7}\addcontentsline
                           {toc}{#1}{\ifnum #2>\c@secnumdepth \else
                             \protect\numberline{\csname the#1\endcsname}\fi
                       #7}}\fi
     \@xsect{#5}}
\def\@begintheorem#1#2{\it \trivlist \item[\hskip \labelsep{\bf #1\ #2.}]}

\def\plain{plain}\ifx\fmtname\plain\csname fi\endcsname
     
     \input docstrip
     \preamble

     Do not distribute the stripped version of this file.
     The checksum in the header refers to the documented version.

     \endpreamble
     \generateFile{here.sty}{t}{\from{here.doc}{}}
     \endinput
\fi
\ifcat a\noexpand @\let\next\relax\else\def\next{%
    \documentstyle[here,doc]{article}\MakePercentIgnore}\fi\next
\ifx\@Hxfloat\@Hundef\else\expandafter\endinput\fi
\let\@Hxfloat\@xfloat
\def\@xfloat#1[{\@ifnextchar{H}{\@HHfloat{#1}[}{\@Hxfloat{#1}[}}
\def\@HHfloat#1[H]{%
\expandafter\let\csname end#1\endcsname\end@Hfloat
\vskip\intextsep\vbox\bgroup\def\@captype{#1}\parindent\z@
\ignorespaces}
\def\end@Hfloat{\egroup\vskip \intextsep}

\makeatother

\catcode`\@=11
\renewcommand{\section}{
        \setcounter{equation}{0}
        \@startsection {section}{1}{\z@}{-3.5ex plus -1ex minus
        -.2ex}{2.3ex plus .2ex}{\large\bf}
        }
\catcode`@=12

\begin{document}

\begin{center}
{\Large 
{\bf  Zero Spacing Distributions for Differenced $L$-Functions}
}

\vspace{1.5\baselineskip}
{\em Jeffrey C. Lagarias} \\
\vspace*{.2\baselineskip}
Dept. of Mathematics \\
University of Michigan \\
Ann Arbor, MI 48109-1109\\
\vspace*{1.5\baselineskip}

(January 1, 2006) \\
\vspace{3\baselineskip}
{\bf ABSTRACT}
\end{center}
The paper studies the local zero spacings of 
deformations of the Riemann $\xi$-function under certain
averaging and differencing operations.
For real  $h$ we consider the entire functions
$A_h(s) := \frac{1}{2}\left(\xi( s+h) + \xi(s-h)\right)$
and 
$B_h(s) = \frac{1}{2i}\left( \xi(s+h) - \xi(s-h)\right).$ 
For  $|h| \ge \frac{1}{2}$ the zeros of $A_h(s)$ and $B_h(s)$
all lie on the critical line $\Re(s) = \frac{1}{2}$
and are simple zeros. The number of zeros 
of these functions to height $T$ has asymptotically 
the same density as the Riemann zeta zeros.
For fixed $|h| \ge \frac{1}{2}$ 
the  distribution of normalized zero spacings 
of these functions up to height $T$ converge as $T \to \infty$
to  a limiting distribution, which consists of  equal 
spacings of size $1$.
That is, these zeros are asymptotically regularly spaced. 
Assuming the Riemann hypothesis, the 
same properties  hold for all nonzero $h$.
In particular,  these
averaging and differencing operations destroy the 
(conjectured) GUE distribution of the zeros of
the $\xi$-function, which should hold at  $h=0$. 
Analogous  results hold for all 
completed Dirichlet $L$-functions $\xi_{\chi}(s)$
having $\chi$ a primitive character.  \\

%
%
%
%
\setlength{\baselineskip}{1.0\baselineskip}

\section{Introduction}

This paper establishes  results on the vertical
spacing distribution of
zeros of certain deformations of Dirichlet $L$- functions
formed using averaging and  differencing  operators. 
Study of these zero spacings is
motivated in part by the GUE conjecture for
vertical zero spacings of $L$-functions, which we
now recall.

There is a great deal of evidence suggesting
that the normalized spacings between the
 nontrivial zeros of the Riemann zeta function
have a ``random'' character described by the
eigenvalue statistics of a random Hermitian matrix
whose size $N \to \infty$. The resulting statistics are
the large $N$ limit of normalized eigenvalue
spacings for random Hermitian matrices drawn from  the GUE distribution
(``Gaussian unitary ensemble''). This limiting distribution is 
identical to the large $N$ limit of  normalized
eigenvalue spacings for random unitary matrices
drawn from the CUE distribution (``circular unitary ensemble''),
i.e. eigenvalues of  matrices drawn from $U(N)$
using Haar measure. (The GUE and CUE spacing  distributions
are not the same for finite $N$.)
More precisely, one compares the normalized
spacings of $k$ consecutive zeros with the 
limiting joint probability distribution
 of the normalized spacings of  $k$ 
adjacent eigenvalues
of random hermitian $N \times N$ matrices,
as $N \to \infty$.
The relation of zeta zeros with random matrix theory was first suggested
by work of H. Montgomery \cite{Mo73} which concerned 
the pair correlation of zeros of the zeta function.
Montgomery's  results showed (conditional on  the Riemnan
hypothesis) that there must be 
some randomness in the spacings of zeros, and were
consistent with the prediction of the GUE distribution.
A. M. Odlyzko \cite{Od87} made  extensive numerical computations
with zeta zeros, now  up to height $T= 10^{22}$, which
show an extremly impressive fit of zeta zero
spacings with predictions of the GUE  distribution.

The GUE distribution  of zero spacings 
is now thought to hold
for all automorphic $L$-functions,
specifically for principal $L$-functions attached
to $GL(n)$, cf.
Katz and Sarnak~\cite{KS99}. Further evidence  for
this was given in 
Rudnick and Sarnak~\cite{RS96}, 
conditionally on a suitable generalized Riemann hypothesis.
They showed that  the evaluation of
consecutive zero gaps against certain test functions (of
limited compact support) agrees with the GUE
predictions. There is also 
supporting numerical evidence for certain 
principal $L$-functions attached to $GL(2)$.

To describe the GUE conjecture, 
we   define the {\em normalized (vertical) zero spacing} 
$\delta_n$ of the $n$-th zero 
$\rho_n=\beta_n + i \gamma_n$ to be
$$
\delta_n := (\gamma_{n+1} - \gamma_n)\frac{1}{2\pi}
\log (\frac{|\gamma_n|}{2 \pi}).
$$
The GUE Conjecture prescribes the limiting distribution of 
any finite set of consecutive normalized zero spacings 
$(\delta_n, \delta_{n+1},...,\delta_{n+k})$, with $k$ fixed
and $n \to \infty$. 
For simplicity we  will consider
the case of consecutive zeros (k=1).

\paragraph{\em GUE Conjecture for Consecutive Zeros.} 
{\em 
The vertical spacings of consecutive normalized zeros 
$\{\delta_n: 1 \le n \le T$
of the Riemann zeta function have a limiting
distribution as $T \to \infty$ given by $p(0,u) du$,
the limiting distribution for consecutive normalized spacings
of eigenvalues  for Gaussian
random $N \times N$  Hermitian matrices as $N \to \infty$.
} \\

\noindent The density $p(0,u)$  is a continuous density
supported on the half-line $u \ge 0$ 
with $p(0,0) = 0$ and which is positive for all $u >0$.
More generally 
there is a density $p(k-1,u)$ giving
the distribution of the spacing to the $k$-th consecutive
normalized zero, which is positive for all $u>0$,
see Mehta and des Cloizeaux \cite{MC72}, Mehta \cite{Me04},
and Odlyzko \cite{Od87}.

At present there is no  known 
structural reason that would explain 
why the normalized zeta zeros might obey the GUE distribution. 
However it is known that 
the  GUE $n$-level spacing distribution is
completely specified by its moments, and
some moments of
the  GUE distribution can be related to (conjectural)
distributions of prime $k$-tuples over various ranges,
see Goldston and Montgomery \cite{GM87} and 
Montgomery and Soundararajan \cite{MS99}.
Recent work of Conrey and Gamburd \cite{CG03} shows
unconditionally that ``pseudo-moments'' of partial sums of the
zeta function on the critical line do exhibit behavior
predicted by the GUE distribution. 

To gain insight into the origin 
and stability of  the GUE property 
it seems useful to study the effect of operators
on functions that preserve 
the property  of having all zeros on a line.
One can then study  the effect of such operators
on the distribution of local  zero spacings.
Here we consider certain sum and difference operators 
constructed using the 
translation operator
${\bf T}_h f(x) := f(x+h)$, such as
${\bf A}_h (f)(x) := \frac{1}{2}\left(f(x+h) + f(x-h)\right)$
and ${\bf B}_h(f)(x):= \frac{1}{2i}\left(f(x+h) - f(x-h)\right)$.
We apply these operators to the Riemann $\xi$-function,
and later to (completed) Dirichlet $L$-functions.

We obtain families of functions by viewing $h$ as a real parameter.
Applied to the  Riemann $\xi$-function
$\xi(s) = \frac{1}{2} s(s-1) \pi^{-\frac{s}{2}}\Gamma(\frac{s}{2})\zeta(s)$, 
we obtain 
the family of ``averaged'' functions
$A_h(s) = \frac{1}{2} \left(\xi(s+h) + \xi(s-h)\right)$,
where $h$ is a real parameter, and another family consists of 
``differenced'' functions
$B_h(s) =  -\frac{1}{2i} \left(\xi(s+h) - \xi(s-h)\right)$.
In fact, one can naturally
insert an second parameter $0 \le \theta < 2 \pi$ and consider
$$
A_{h, \theta}(s) :=  \frac{1}{2} \left(\cos \theta (\xi(s+h) +\xi(s-h))
+i \sin \theta (\xi(s+h)- \xi(s-h)) \right).
$$
There is  an analogous extended family  $B_{h, \theta}(s)$
associated to  $B_h(s)$.

In \S2 we  show 
unconditionally that  
for fixed  $|h| \ge \frac{1}{2}$
the zeros of $A_h(s)$ and $B_h(s)$ lie 
on the critical line $\Re(s) = \frac{1}{2}$
and are simple zeros, and they interlace.
 Assuming the Riemann hypothesis, we show that 
the same result holds for all nonzero  $h$.

In \S3 we show that for all real $h$ that
the functions in this family have 
the same asymptotic density of
zeros as the zeta function. The results together
permit a definition of normalized zero spacings for the
functions in this family.

In \S4 we establish the main results of the paper.
We show that for non-zero $h$ the functions 
$A_{h, \theta}(s)$ and $B_{h, \theta}(s)$ have a limiting
normalized zero spacing distribution for any fixed
number $k$ of consecutive zeros,  which corresponds
to a delta measure with equal spacings of size $1$. 
This  is proved unconditionally when $|h| \ge \frac{1}{2}$,
and proved conditionally on  the Riemann hypothesis for the remaining
range $0 < |h| < \frac{1}{2}.$
These results assert that
for $|h| \ne 0$  the normalized zero spacings are completely regular.
We conclude that these differencing operations destroy the 
GUE property. Of course for $h=0$ the GUE distribution is
expected to hold.

In \S5 we observe that these  results also 
hold for zeros of  (completed) Dirichlet $L$-functions $\xi(s, \chi)$
for primitive characters $\chi$,
with little change in the proofs.
The proof method should also extend to various automorphic 
$L$-functions for  $GL(N)$, for $N \ge 2$.

In \S6 we interpret the results of \S2  in
terms of the de Branges theory of Hilbert spaces of
entire functions, which  was another motivation
for this work. The notations used in this paper were
chosen to be compatible with de Branges's theory,
under the change of variable $s= \frac{1}{2} - iz$.
The interpretation is  that these sum
and differencing operators applied to  L-functions produce
structure functions of de Branges Hilbert spaces
of entire functions 
when $|h| \ge \frac{1}{2}$,
and, conditionally on the Riemann hypothesis, for all nonzero $h$.
For example, Lemma~\ref{le21} asserts 
that $E_h(\frac{1}{2}-iz):= \xi(\frac{1}{2}+ h - iz)$ 
is a de Branges structure function, under these hypotheses.
The particular case with $h=\frac{1}{2}$
gives the structure function $E(z) = \xi(1-iz)$
which was discussed  
in 1986 by de Branges \cite{deB86}.
We describe some features of the de Branges theory, and observe
that it provides a ``spectral'' interpretation of
the zeros $\frac{1}{2} + i \gamma$ 
of the functions
$A_{h, \theta}(s)$
(resp. $B_{h, \theta}(s)$).
It produces for each $\theta$ 
an unbounded self-adjoint operator on a Hilbert
space having the  negative
imaginary parts of the zeros (that is,  $-\gamma$) as its 
spectrum. 

In \S7 we make concluding remarks and suggest some 
directions for further work.  In particular 
there likely exist  random matrix analogues of these 
results.

We add some comments regarding the proofs.
In Theorem~\ref{th21} we
establish simplicity of zeros, all lying on the critical line,
using a function-theoretic 
result of de Branges~\cite{deB59} which 
applies generally to all de Branges structure functions.
For completeness we  include a self-contained proof of 
this result (Lemma~~\ref{le22}).
The results of \S3 are standard. The main result of \S4 is
proved by viewing the zeros of 
$E_{h, \theta}(s)= A_h(s) - i B_h(s)$ as perturbations
of zeros of the real part of the archimedean factor contribution
on the critical line.
These zeros are regularly spaced, and  under
the given hypotheses  the perturbations are shown to be sufficiently small
to not significally disturb these spacings. 

Concerning prior work, Xian-Jin Li informs me that the results 
in \S2 and \S5  were known to  de Branges
in the late 1980's. These include Lemma 2.1, Theorem 2.1 and
their extensions to Dirichlet $L$-functions stated in \S5.
De Branges reportedly
covered some of these results 
in a lecture course on Hilbert spaces of entire functions 
given at Purdue in 1987-1988. In this regard,
the  present paper serves to make proofs of these results 
available. Recently Haseo Ki ~\cite{Ki04} obtained 
results analogous to those in \S2 for differencings of the
meromorphic function 
$\hat{\zeta}(s)= \pi^{-\frac{s}{2}} \Gamma(\frac{s}{2} )\zeta(s)
= 2 \frac{\xi(s)}{s(s-1)}$,
e.g. for $h \ge \frac{1}{2}$ all zeros of
$\tilde{A}_h(s)= \frac{1}{2}\left( \hat{\zeta}(s+h)+ 
\hat{\zeta}(s-h)\right)$ 
lie on the critical line. Ki's methods should also
apply to the function $\xi(s)$.

\paragraph{Acknowledgments.} Much of this work was
done when the author worked at AT\&T Labs,
and he thanks AT\&T Labs for research support. 
He thanks  D. Cardon, S. Gonek, H. Ki, 
M. Rubinstein, K. Soundararajan  and the 
reviewer for helpful comments and references.
In particular M. Rubinstein 
pointed out the relevance of a result
in Titchmarsh cited in Lemma~\ref{le41}(3) below.

%
%
%
%
\section{Differenced Riemann $\xi$-function}

The Riemann $\xi$-function $\xi(s)$ is given by
$
\xi(s) = \frac{1}{2} s (s-1) \pi^{- \frac{s}{2}} 
\Gamma(\frac{s}{2}) \zeta(s).
$
It is real on the real axis. More importantly here, it
is real on the critical line $\Re(s) = \frac{1}{2}$,
and satisfies the functional equation $\xi(s) = \xi(1-s)$.
Its zeros are confined to the (open) critical strip
$0 < \Re(s) < 1$, and as an
entire function of order $1$ it has the Hadamard factorization
$$
\xi(s) = e^{A(\chi_0) + B(\chi_0)s }
\prod_{\rho} \left( 1 - \frac{s}{\rho}\right) e^{\frac{s}{\rho}}.
$$
The  spacing of its zeros is asymptotically regular, and this
permits one to
show that the convergence factors $e^{s/\rho}$ can be
factored out into the lead term,  if the product is taken in 
a suitable order,
for example considering  partial products over
all terms with $|\rho| < T$ and taking a limit as $T \to \infty$.
In that case one obtains a modified Hadamard
product representation
\beql{200}
\xi(s) = e^{A^{\ast}(\chi_0) + B^{\ast}(\chi_0)s}
\prod_{\rho, \ast} \left( 1 - \frac{s}{\rho}\right)
\eeq
where the asterisk in the product indicates that its terms
are to be combined in a suitable order.
For the $\xi$ function it suffices to group all
zeros in pairs $(\rho, 1-\rho)$ if on the critical line or in 
quadruples as $(\rho, 1 - \rho, \bar{\rho}, 1 - \bar{\rho})$ if not 
to get absolute convergence of the 
modified Hadamard product \eqn{200}. We have
$$
\xi(s) = \lim_{T \to \infty} e^{A^{\ast}(\chi_0) + B^{\ast}(\chi_0)s}
\prod_{|\rho| < T} \left( 1 - \frac{s}{\rho}\right).
$$
It is known that 
$$
A(\chi_0)= A^{\ast}(\chi_0) = \log \xi(0) = - \log 2
$$
and $B(\chi_0)$ is real with 
$$
B(\chi_0) = - \frac{1}{2} \gamma -1 + \frac{1}{2} \log 4 \pi \approx -0.023,
$$
where $\gamma$ is Euler's constant,
see Davenport~\cite[p. 83]{Da80}. In addition
$B(\chi_0) = - \sum_{\rho} \frac{1}{\rho}$ where the sum
takes the zeros in complex conjugate pairs, and from this one
deduces that 
\begin{equation}~\label{200a}
B^{\ast}(\chi_0) = 0.
\end{equation}
For the purposes here what is important is that $\Re(B^{\ast}(\chi_0))\ge 0$,
see Lemma~\ref{le21}.

In this section we consider for a real parameter $h$ the functions
$$
A_h(s) := \frac{1}{2}\left(\xi(s+h) + \xi(s-h)\right)
$$
and
$$
B_h(s) :=  -\frac{1}{2i}\left( \xi(s+h) - \xi(s-h)\right).
$$
We have
$$ 
\xi(s-h) = \overline{\xi(s+h)},
$$
which follows from the reflection symmetry of
$\xi(s)$ around the line $\Re(s) = \frac{1}{2}$, i.e. 
$$
 \xi (\frac{1}{2} -x + iy) =\overline{\xi(\frac{1}{2} +x -iy)}.
$$
It follows that $A_h(s)$ and $B_h(s)$ are both real
on the critical line and satisfy there
$$
A_h(\frac{1}{2} + it) = \Re \left( \xi(\frac{1}{2} + h + it)\right)
$$
and
$$
B_h(\frac{1}{2} + it) = -\Im\left(\xi(\frac{1}{2} + h + it)\right).
$$
In particular $A_0(s) = \xi (s)$.
It is also immediate that $A_{-h}(s) = A_h(s)$
and $B_{-h}(s) = - B_{h}(s)$, so that without
loss of generality we need only consider $ h \ge 0$.
The fact that $\xi(s)$ is real on the real axis gives rise to the
extra symmetry $A_h(\bar{s}) = \overline{A_h(s)}$,
and similarly for $B_h(s)$.
%
%
\begin{lemma}~\label{le21}
(1) If $h \ge \frac{1}{2}$, then
\beql{201}
|\xi(h+s)| > |\xi(h+1-\bar{s})| ~~~\mbox{for}~~ \Re(s) > \frac{1}{2}.
\eeq
(2) Assuming the Riemann hypothesis, 
the inequality \eqn{201} holds for each $h > 0$.
\end{lemma}

\paragraph{Remark.} The inequality \eqn{201} can be 
stated alternatively (and apparantly more elegantly) as \\
$$|\xi(h+s)| > |\xi(h + (1-s))|,$$
making use of the reflection symmetry of $\xi(s)$ around
the real axis, i.e. $\xi(\bar{s}) = \overline{\xi(s)}$.
However the restated  form does not generalize to
Dirichlet $L$-functions, which do not all have this symmetry.
In the form stated the lemma generalizes to Dirichlet
$L$-functions, with the relevant property being that
the (completed) $L$-function have constant phase on
the critical line.

\paragraph{Proof.} 
We show that the inequality holds 
term by term for each factor in the  modified Hadamard product \eqn{200}
for $\xi(s)$. For the exponential factor, let $s = \sigma + it$,
and we have
$$
|e^{B^{\ast}(\chi_0)(h+s)}|= e^{\Re(B^{\ast}(\chi_0))(h + \sigma)}
e^{- \Im(B^{\ast}(\chi_0))t}
$$
$$
|e^{B^{\ast}(\chi_0)(h+1-\bar{s}}|= e^{\Re(B^{\ast}(\chi_0))(h +1- \sigma)}
e^{- \Im(B^{\ast}(\chi_0))t}.
$$
The condition $\Re\left( B^{\ast}(\chi_0)\right) \ge 0$ is sufficient to
imply that for 
$\sigma > \frac{1}{2}$, 
$$
|e^{B^{\ast}(\chi_0)(h+s)}| \ge |e^{B^{\ast}(\chi_0)(h+1-\bar{s}}|,
$$
as desired. In fact  $ B^{\ast}(\chi_0)=0$ by \eqn{200a}.
 
Now we consider the product factors for each zero separately.
Let $\rho= \beta + i \gamma$ be a zero
of $\xi(s)$, so that $0 < \beta < 1$, and under
the Riemann hypothesis $\beta = \frac{1}{2}$.
Now set $s= \sigma + it,$ where we will
suppose $\sigma > \frac{1}{2}$, and we will show under
the stated hypotheses that
\beql{203}
|1 - \frac{h+s}{\rho}| > |1 - \frac{h+1 - \bar{s}}{\rho}|.
\eeq
Multiplying by $|\rho|$ gives the equivalent inequality
$$
|\beta - h - \sigma + i(\gamma - t)| > 
|\beta - h -1 + \sigma + i(\gamma - t)|,
$$
and this in turn is equivalent to the inequality
\beql{204}
|\beta -h - \sigma| > |\beta - h -(1 - \sigma)|.
\eeq

(1) Suppose $h \ge \frac{1}{2}$ so that
$\beta -h -\frac{1}{2} <0$. Then one has
$$
|\beta - h - \sigma| = |(\beta -h - \frac{1}{2}) - (\sigma - \frac{1}{2})|=
|\beta - h - \frac{1}{2}| + |\sigma - \frac{1}{2}|,
$$
while
$$
|\beta -h -1 + \sigma| = |(\beta - h - \frac{1}{2}) + (\sigma - \frac{1}{2})|
< |\beta - h - \frac{1}{2} | + |\sigma - \frac{1}{2}|.
$$
The  application of the triangle inequality is strict because both
terms in it are nonzero and have opposite signs. This establishes
\eqn{204}, hence \eqn{203}, in this case.

(2) Now assume RH, so that $\beta = \frac{1}{2}$, and suppose $h >0$. Now 
$\beta - \sigma = \frac{1}{2} - \sigma < 0$,
and then
$$
|\beta - h - \sigma| = | (\beta - \sigma) -h| = |\beta - \sigma| + |h|,
$$
while
$$
|\beta - h -(1 - \sigma)| = |(\sigma - \frac{1}{2}) - h|
< |\beta - \sigma| + |h|.
$$
The application of the triangle inequality in the last line is
strict because  both terms in it are nonzero and of opposite
sign. This establishes \eqn{204}, hence \eqn{203}, in this case.
 ~~~~$\bsq$. \\

The next lemma, due to de Branges ~\cite[Lemma 5]{deB59}, 
formulates  a basic property underlying the de Branges
theory of Hilbert spaces of entire functions. The de Branges
theory is formulated in terms of a variable $z$ with $s= \frac{1}{2} -iz$
and considers functions satisfying $|E(z)| > |E^{\sharp}(z)|$
when $\Im(z) > 0$, where $E^{\sharp}(z) := \overline{E(\bar{z})}$
is an involution acting on entire functions.  
Here we re-express de Branges's result in terms of the $s$-variable, and 
the involution becomes
$E^{\sharp}(s) = \overline{E( 1 - \bar{s})}$.  

%
\begin{lemma}~\label{le22}
Let $E(s)$ be an entire function that satisfies 
\beql{206}
|E(s)| > |E(1 - \bar{s})|~~~\mbox{when}~~ \Re(s) > \frac{1}{2}.
\eeq
Write 
$E(s) = A(s) - i B(s)$ with
$A(s) = \frac{1}{2} (E(s) + \overline{E(1- \bar{s})})$
and $B(s) = -\frac{1}{2i} (E(s) - \overline{E(1- \bar{s})})$,
so that $A(s)$ and $B(s)$ are real-valued on the critical
line $\Re(s) = \frac{1}{2}.$
Then $A(s)$ and $B(s)$ have all their zeros lying on 
the critical line $\Re(s) = \frac{1}{2}$, and these
zeros interlace.
\end{lemma}

\paragraph{Remarks.} (1) The interlacing property allows
$A(s)$ and $B(s)$ to have multiple zeros. 
By ``interlacing of zeros'' we mean there
is a numbering of zeros of the two functions, 
$\{\rho_n(A)= \frac{1}{2} + \gamma_n(A) : n \in \ZZ \},$ resp. 
$\{\rho_n(B)= \frac{1}{2} + \gamma_n(B) : n \in \ZZ \},$
with imaginary parts in  increasing
order,  counting
zeros with multiplicity, such that 
$$
\gamma_n(A) \le \gamma_n(B) \le \gamma_{n+1}(A)
$$
holds for all allowed $n$.

(2) There are no growth restrictions on the maximum modulus of
functions $E(s)$ in Lemma~\ref{le22}. It 
can be shown that there exist
entire functions $E(s)$ of fast growth, for example
of infinite order, satisfying  the
condition \eqn{206}.

\paragraph{Proof.}
The definition of $A(s)$ and $B(s)$ shows that 
for real $\alpha$ they have the symmetry
\beql{206a}
A(\frac{1}{2} - \alpha + it) =
\overline{A(\frac{1}{2} + \alpha + it)},
\eeq
and similarly 
$B(\frac{1}{2} - \alpha + it) =
\overline{B(\frac{1}{2} + \alpha + it)}.$
Consequently $A(s)$ and $B(s)$ are real-valued on the
critical line $\Re(s) = \frac{1}{2}.$

To see that $A(s)$ has all its zeros on the critical line,
we observe that for $\Re(s) > \frac{1}{2}$, \eqn{206} gives
$
|A(s)| = |E(s) + \overline{E(1 - \bar{s})}|
\ge |E(s)| - |E(1 - \bar{s}| > 0.$
If $\Re(s) < \frac{1}{2}$ then $\Re(1 - \bar{s}) > \frac{1}{2}$
and \eqn{206} gives
$$
|A(s)| = |E(s) + \overline{E(1 - \bar{s})}|
\ge |E(1- \bar{s})| - |E(1 - \overline{(1-\bar{s})}| > 0.
$$
This establishes that  all zeros of $A(s)$ lie on the critical line,
and the proof that $B(s)$ has the same property is similar.

To see that the zeros of $A(s)$ interlace with those of
$B(s)$ on the critical line (counting multiplicities),
we observe that the  property \eqn{206}  implies for real $\alpha > 0$ that
$$
\lim_{\alpha \to 0^{+}}
\frac{ |E( \frac{1}{2} + \alpha + it)|^2 - 
|E( \frac{1}{2} - \alpha+ it)|^2}{ 2\alpha} \ge 0.
$$
On inserting $E(s)$ and its conjugate $\overline{E(s)}$
into this and using the symmetry \eqn{206a}
for $A(s)$ and $B(s)$
a calculation yields, for real $\alpha > 0$, on letting
$s_{\alpha} = \frac{1}{2} + \alpha +it$, 
\begin{eqnarray*}
\frac{ |E( \frac{1}{2} + \alpha + it)|^2 - 
|E( \frac{1}{2} - \alpha+ it)|^2}{ 2\alpha} &=&
i \left( \frac{ A(s_{\alpha}) \overline{B(s_{\alpha})} - 
\overline{A(s_{\alpha})} B(s_{\alpha})}{\alpha}
\right) \\
&=& 
-i A(s_{\alpha}) 
\left(\frac{B(s_{\alpha}) - \overline{B(s_{\alpha})} } {\alpha}\right)  \\
&&  +
i B(s_{\alpha})
\left(\frac{A(s_{\alpha}) - \overline{A(s_{\alpha})}} {\alpha}\right). 
\end{eqnarray*}
Letting $\alpha \to 0^{+}$ we  deduce
that for $s = \frac{1}{2} + it$, there holds 
\beql{207}
(i\frac{d}{ds} A(s)) B(s) - (i\frac{d}{ds} B(s)) A(s) \ge 0.
\eeq
The left hand side of this inequality is 
real-analytic in $t$ so it can vanish only 
at isolated points, provided  it is not identically zero,
and this would give the interlacing
property. It remains to establish that it is not identically zero.

We can rephrase the result \eqn{207} in terms of a 
continuous real-valued ``phase function'' $\phi(t)$  defined by 
\beql{207a}
E(s) := |E (\frac{1}{2} + it)| e^{i \phi(t)},
\eeq
where we choose  the
unique allowed value $0 \le \phi(0) < 2\pi$, and
then extend it to all $t$ by continuity.
Now \eqn{206} implies that  $E(s)$ has no zeros in $\Re(s) > \frac{1}{2}$
so we can define a single-valued version of $\log E(s)$ there,
such that $\phi(t)= \Im(\log  E(\frac{1}{2} + it))$. Using
$A(s) = |E(s)| \cos \phi(t)$ and $B(s) = - |E(s)| \sin \phi(t)$
and $\frac{d}{ds} = i \frac{d}{dt}$ on the line $s= \frac{1}{2} +it$
one finds that
the inequality \eqn{207} simplifies to 
\beql{208}
|E(\frac{1}{2} +it)|^2 \frac {d}{dt}  \phi(t) \ge 0.
\eeq
To show the left side is not identically zero, 
it suffices to  show that
$\phi(t)$ is nonconstant. We argue by contradiction.
If it were constant, say
$E(s) = |E(s)|e^{i\alpha}$, then the function 
$\tilde{E}(s)= E(s) e^{-i\alpha}$  would be
real-valued on the critical line, 
and the reflection principle would then give 
$\tilde{E}(1-\bar{s})=\overline{\tilde{E}(s)}$, which contradicts
\eqn{206}. We conclude that the real-analytic
function $\phi(t)$ is nonconstant, so
$\frac {d}{dt}  \phi(t)$
 has isolated zeros, and \eqn{208} implies 
 that the function $\phi(t)$ is a strictly increasing
function of $t$. This certifies that 
the zeros of $A(t)$ and $B(t)$ interlace.
~~~$\bsq$ \\

If
 a  function $E(s)= A(s) - i B(s)$ satisfies the hypothesis of 
Lemma~\ref{le22}
then so does the function $E_{\theta}(s) = e^{i \theta}E(s)$, when
$0 \le \theta < 2 \pi.$ It has the decomposition
$E_{\theta}(s) = A_{\theta}(s) - i B_{\theta}(s)$, where
$$
A_{\theta}(s) =  ( \cos \theta) A(s) + (\sin \theta) B(s) 
$$
$$
B_{\theta}(s) = -(\sin \theta) A(s) + (\cos \theta) B(s).
$$
In particular,  $B_{\theta}(s) = A_{\theta + \pi/2}(s).$

We now specialize to the case of the shifted $\xi$-function
$E_h(s) = \xi(s+h)$, for which  $E_h(s) = A_h(s) - i B_h(s)$
with 
$$
A_h(s) := \frac{1}{2}\left( \xi(s+h) + \overline{\xi(1-\bar{s}+h)}\right)= 
\frac{1}{2}\left(\xi(s+h)+\xi(s-h)\right)
$$
$$
B_h(s) := -\frac{1}{2i}\left(\xi(s+h) - \overline{\xi(1-\bar{s}+h)}\right)= 
-\frac{1}{2i}\left(\xi(s+h)-\xi(s-h)\right).
$$
Here we used  $\xi(\bar{s})= \overline{\xi(s)}$ and
the functional equation $\xi(s)= \xi(1-s)$.
We consider  
$$E_{h, \theta}(s) = e^{i \theta}E_h(s) = e^{i \theta} \xi(s+h),$$
and obtain
\begin{eqnarray*}
A_{h, \theta}(s) & := & (\cos \theta) A_h(s) +( \sin \theta) B_h(s) \\
&=& \frac{1}{2}(\cos \theta) \left(\xi(s+h) +\xi(s-h)\right)
-\frac{1}{2i} (\sin \theta) \left(\xi(s+h)- \xi(s-h)\right) \,
\end{eqnarray*}
\begin{eqnarray*}
B_{h, \theta}(s) & := & (-\sin \theta) A_h(s) + (\cos \theta) B_h(s) \\
&=&-\frac{1}{2}(\sin \theta) \left(\xi(s+h) +\xi(s-h)\right)
 -\frac{1}{2i}(\cos \theta) \left(\xi(s+h)- \xi(s-h)\right) .
\end{eqnarray*}
For $\theta=0$  we recover 
$A_{h,0}(s) = A_h(s)$ and 
and  $B_{h,0}(s) = B_h(s)$.


\begin{theorem}~\label{th21}
(1) For $|h| \ge \frac{1}{2}$ and any $0 \le \theta < 2 \pi$,
the entire functions
$A_{h, \theta}(s)$ and $B_{h, \theta}(s)$ have all their
zeros on the critical line $\Re(s) = \frac{1}{2}.$
These zeros are all simple zeros, and they interlace.

(2) Assuming the Riemann hypothesis, for
$0 < |h| < \frac{1}{2}$ and any $0 \le \theta < 2 \pi$
the functions
$A_{h, \theta}(s)$ and $B_{h, \theta }(s)$ have 
all their zeros on the critical line $\Re(s) = \frac{1}{2}.$ 
These zeros are all simple zeros, and they interlace.
\end{theorem}

\paragraph{Proof.}

Lemma~\ref{le21} shows that under the stated assumptions
the function
$E(s) := E_{h}(s) = \xi(s+h)$ 
satisfies the hypothesis \eqn{206}
of Lemma~\ref{le22}, with $A(s) = A_{h}(s)$ and 
$B(s) = B_{h}(s).$
The same holds 
more generally for $E_{h, \theta} (s)= e^{i \theta} E_h(s)$.
Applying Lemma~\ref{le22} we  conclude in both cases (1) and (2) that 
$A_{h, \theta} (s)$
and $B_{h, \theta}(s)$ have all their zeros on the critical line
$\Re(s) = \frac{1}{2}$, and  these zeros interlace.

To conclude that these zeros are all simple zeros, 
we argue by contradiction. Suppose that one of $A_{h, \theta}(s)$ or
$B_{h, \theta}(s)$ had a multiple zero at a point $s_0$.  
Then by the interlacing property  they necessarily
have a common zero at $s_0$. Since this point is on the
critical line, we infer that the real and imaginary parts
of $e^{i \theta}\xi(s_0+h)$ both vanish, whence $\rho_0 = s_0+h $
is a zero of $\xi(s)$. In case (1), 
$\Re(\rho_0) = \Re(s_0) + h= \frac{1}{2} +h \ge 1$,
a contradiction. In case (2), assuming RH, $h >0$
gives $\rho_0 = \frac{1}{2} + h > \frac{1}{2}$, a
contradiction.
~~~$\bsq$ \\

%
%
%
%

\section{Global Asymptotics of Zeros}

We determine the 
asymptotic zero density of functions $A_{h, \theta}(s)$ and 
$B_{h, \theta}(s)$.
The proof method is standard, using the argument principle.
 Let $N(T, F)$ count the
number of zeros of $F(s)$ having $|\Im(s)| \le T$.
%
%
\begin{theorem}~\label{th31}
(1) For $|h| \ge \frac{1}{2}$ and any $0 \le \theta < 2 \pi$,
then for all  $|T| \ge 2$ there holds
\beql{301}
N(T, A_{h, \theta}) = \frac{1}{\pi} T \log T - 
\frac{1}{\pi}(\log (2\pi) + 1)T + 
O( \log T),
\eeq
where the implied constant in the O-notation is independent of $h$
and $\theta$.
A similar formula holds for $N(T, B_{h, \theta}).$

(2) Assuming the Riemann hypothesis, the formula \eqn{301}
for $N(T, A_{h, \theta})$ is valid for all nonzero $h$.
A similar formula holds for $N(T,B_{h, \theta})$
for all nonzero $h$. 
\end{theorem}

\paragraph{Proof.}
We study $\xi(s+h)$ in a half-plane $\Re(s+h) > \sigma_0$ in which
it has no zeros, where the value of $\sigma_0 =1$ unconditionally
and $\sigma_0= \frac{1}{2}$ if the Riemann hypothesis is
assumed. Under these conditions we can choose
a single-valued branch of  $\log \xi(s)$ in this half-plane,
and we choose the branch which is real-valued on the real axis.
We set
$$
\xi(\frac{1}{2} + h + it) = |\xi(\frac{1}{2} +h +it)| e^{i \varphi_h(t)}
$$
in which
$$
\varphi(t) = \arg( \xi(\frac{1}{2} + h + it)) = 
\Im \left(\log \xi(\frac{1}{2} + h + it) \right).
$$
The function $E_h(s)=\xi(s+h)$ satisfies the conditions of
Lemma~\ref{le22} and its proof shows that $\varphi(t)$ is a
strictly increasing function of $t$. By Theorem~\ref{th21} the
zeros of $A_h(s)$ and $B_h(s)$ are simple and lie on the
critical line. We write the zeros of $A_{h, \theta}(s)$  as 
$\rho_n(A_{h, \theta}) = \frac{1}{2} + i \gamma_n(A_{h, \theta})$,
 enumerated in
order as 
$$ ...< \gamma_{-1}(A_{h,\theta}) < 0 \le \gamma_0(A_{h,\theta}) < 
\gamma_1(A_{h,\theta}) < ...,$$
and we denote the zeros of $B_{h,\theta}(s)$ similarly,
as $\rho_n(B_{h,\theta}) = \frac{1}{2} + i \gamma_n(B_{h,\theta}).$
In terms of the argument $\varphi_h(t)$ we 
 have that these zeros comprise all the solutions to  
$$
\varphi_h\left(\gamma_n(A_{h,\theta}) \right) \equiv \frac{\pi}{2}+ \theta   
~(\bmod~\pi)
$$
and
$$
\varphi_h \left(\gamma_n (B_{h, \theta}) \right) \equiv \theta ~(\bmod~\pi).
$$
In particular $\varphi_h \left(\gamma_n(A_{h,\theta})\right) = n\pi + O(1)$
and similarly for $\varphi_h \left(\gamma_n(B_{h, \theta})\right)$.

Because all zeros are on the critical line and the
argument $\varphi(t)$ is strictly increasing there,
to estimate $N_h(T)$ it suffices to bound the change in
argument. By definition  we have  
\beql{406}
\arg(\xi(s)) = \arg(s(s-1)) +  \arg(\pi^{-\frac{s}{2}} \Gamma(\frac{s}{2})) 
+ \arg(\zeta(s)).
\eeq
The contribution to the argument mainly
comes from the Gamma function factor, which we deal with first.

We let $\tilde{\varphi}_h(t)$ denote the argument of
$G_h(t) = \pi^{- (\frac{1}{4} + \frac{h}{2} + \frac{it}{2})} 
\Gamma(\frac{1}{4} + \frac{h}{2} + \frac{it}{2}),$
normalized by the condition that it be zero on the positive real 
axis.
It can be checked that $\tilde{\varphi}_h(t)$ is a 
strictly increasing function of $t$.
Stirling's formula is valid in any sector 
$-\pi + \delta < arg(s) < \pi - \delta$ with $\delta > 0$,
 and asserts that
\beql{403}
\log \Gamma(s) = (s - \frac{1}{2}) \log s - s + \frac{1}{2} \log 2\pi
+ O (\frac{1}{|s|}),
\eeq
where the $O$-constants depend on $\delta > 0$. We choose 
$\delta= \frac{\pi}{2}$, so can  use the formula 
on the half-plane $\Re(s) > 0$. Letting $s = \frac{1}{2}+ h + it$,
and using $\arg(\pi^{-\frac{s}{2}}) = -(\frac{1}{2} \log 2\pi) t + O(1)$
we obtain, for $ t \ge 0$, that 
\beql{404}
\tilde{\varphi}_h (t)= \frac{t}{2} \log \frac{t}{2} 
-\frac{t}{2} (\log 2\pi + 1) + 
\frac{\pi}{2}\left( -\frac{1}{4} + \frac{h}{2} \right) + O (\frac{1}{|t|+1}),
\eeq
and we have $\tilde{\varphi}_h (-t)= -\tilde{\varphi}_h (t).$

We estimate the remaining terms on the right side of \eqn{406}
The first term on the right for  $t \ge 1$ satisfies
\beql{407}
\arg(s(s-1)) = \pi + O\left(\frac{1}{|t|+1}\right).
\eeq
On taking $s= \frac{1}{2} + h + it$, we deduce that
$$
\varphi_h(t)  = \tilde{\varphi}_h(t) + \pi + 
\arg(\zeta(\frac{1}{2} + h + it))
+ O (\frac{1}{|t| + 1}).
$$
It is a standard estimate that in a half-plane $\Re(s) > \frac{1}{2} + h$
containing no zeros there holds
$$
\arg(\zeta(\frac{1}{2} + h + it)) = O ( \log |t|),
$$
see Titchmarsh~\cite[Theorem 9.4]{TH86} on Davenport~\cite[Chap. 15]{Da80}.
(It is based on bounding $\Re\left(\frac{\zeta'}{\zeta}(s) \right).$)
Substituting these estimates in \eqn{406} together with \eqn{404}
yields, for $|t| \ge 2$, that 
$$
\varphi_h(t) = \frac{t}{2} \log \frac{t}{2} -
\frac{t}{2} (\log 2\pi + 1) + O ( \log |t| ).
$$
Using 
$$
|N_h(T) - \frac{1}{\pi}| (\varphi_h(T) - \varphi_h(-T))| \le 2,
$$
the bounds (1) and (2) follow.
~~$\bsq$ \\

We remark that for $0 < h < \frac{1}{2}$
the asymptotics of zeros of 
the real part (resp. imaginary part) of
$\xi(s)$ can be established unconditionally,
with a worse error term. 
Haseo Ki ~\cite{Ki04}
showed analogous results for the meromorphic function 
$\hat{\zeta}(s)= 2\frac{\xi(s)}{s(s-1)}$; an earlier
result of Levinson \cite{Le71}
obtained the lower bound.
Let $\tilde{N}_h(T)$ count the number of zeros of the difference
function 
$$
\tilde{F}_h(s) := \frac{1}{2}\left(\hat{\zeta}(s+h)+ \hat{\zeta}(s-h)\right)
$$
having  imaginary part between $0$ and $T$. Then 
Ki \cite[p. 290]{Ki04} shows that  for $h >0$ 
\begin{equation}~\label{341}
\tilde{N}_h(T) = \frac{T}{2\pi} \log \frac{T}{2\pi} - \frac{T}{2\pi}
+O( \tilde{R}_h(T) \log T),
\end{equation}
in which $\tilde{R}_h(T)$ counts the number of zeros of $\zeta(s)$
in $\Re(s) > \frac{1}{2}+h$ with imaginary part between $0$ and $T$.
It is known that $\tilde{R}_h(T)= 
O \left( T^{3\frac{1-2h}{3-2h}}(\log T)^5\right)$
for $0< h < \frac{1}{2}$
(Titchmarsh \cite[Theorem 9.19(B)]{TH86}). Ki's  methods 
apply to  give similar results for 
the function $\xi(s)$ as well, which unconditionally imply \eqn{301}
replaced with the larger remainder term  $O ( \tilde{R}_h(T) \log T)$.

%
%
%
%

\section{Local Spacing Statistics of Zeros}

In general we can  define local spacing statistics of
zeros of an entire function $G(z)$ having
all its zeros on the real axis, 
as follows. Let $\{\gamma_i : i \in  \ZZ \}$
enumerate the zeros in increasing order, with $\gamma_{0}$
denoting the zero nearest the origin, and let
$N(T)$ denote the number of such zeros
in $-T \le x \le T$. We consider  the
set of  differences of consecutive zeros 
$$
\Sigma(T) :=\{ (\gamma_{n+1} - 
\gamma_{n})\frac{N(T)}{T}: | \gamma_n| \le T\},
$$
normalized to have average spacing $1$
over this interval.  We put
the uniform  probability distribution to $\Sigma_T$, i.e.
assign  equal weights to each observation. We say 
that the function $G(z)$ has {\em limiting local spacing statistics}
if as $T \to \infty$ the distributions $\Sigma_T$
weakly converge to a limiting probability distribution $\Sigma$
on the real line. 

We can apply this definition of local spacing
distributions to the functions here by using the linear
change of variable $s= \frac{1}{2} - iz$ which moves the critical
line to the real axis. 
The definition above,
which works for an arbitary
set of zeros,  differs slightly from  the normalization of  zero spacings
used in Odlyzko~\cite{Od87}. Recall that  he used the 
{\em normalized consecutive zero spacing}
$$
\delta_n := (\gamma_{n+1} - \gamma_n)\frac{1}{2\pi} 
\log \frac{|\gamma_n|}{2\pi}.
$$
Using the asymptotics
  $N_h(T) = \frac{1}{\pi} T \log \frac{T}{2\pi e} + O(\log T)$
given by Theorem ~\ref{th31}, we have
$\gamma_n = 2 \pi \frac{n}{\log~n} + O( \frac{n}{(\log~n)^2})$
as $n \to \infty$, and  this yields 
$$ 
\delta_n =  (\gamma_{n+1} - \gamma_{n})\frac{N(T)}{T} + 
O(\frac{1}{\log\log T}),
$$
valid over the range $T \le n \le N(T)$.
This bound is sufficient to imply that the set
$$
\Sigma^{\ast}(T) : = \{ \delta_n :~ |\gamma_n| < T \}
$$
assigned the uniform probability distribution
on its $N(T)$ elements,
will have the same limiting distribution as that of $\Sigma(T)$ as 
$T \to \infty$. This holds
 in the sense that, if either $\Sigma(T)$ or $\Sigma^{\ast}(T)$
 has a limiting distribution as $T \to \infty$, 
then so does the other, and they agree. 

The main observation of this paper is the
following.

\begin{theorem}~\label{th41} Let $k \ge 1$ be given.

(1) For $|h| \ge \frac{1}{2}$ and $0 \le \theta < 2 \pi$ the functions
$A_{h, \theta}(s)$ and $B_{h, \theta}(s)$ both have a 
limiting joint distribution as $T \to \infty$ of their set of $k$ 
normalized consecutive zero spacings
$(\delta_n, \delta_{n+1}, ..., \delta_{n+k})$. This distribution
of normalized zero spacings consists  of a delta function
located at $(x_1,x_2, ..., x_{k})=(1,1,\cdots, 1)$, 
the ``trivial'' distribution.

(2) Assuming the Riemann hypothesis, the same
results also hold for $0 < |h| < \frac{1}{2}$
and $0 \le \theta < 2 \pi$, 
with the limiting distribution of normalized zero 
spacings existing and being
the ``trivial'' distribution.
\end{theorem}

The idea of the proof is that the zeros  for
$\xi_h(\frac{1}{2} + it)$ can be viewed as a 
perturbation of ``zeros''  associated to the
argument of the $\Gamma$-factor term
$$G_h(t) 
:= \pi^{- (\frac{1}{4} + \frac{h}{2} + \frac{it}{2})} 
\Gamma(\frac{1}{4} + \frac{h}{2} + \frac{it}{2}).
$$
The ``zeros'' are values  $arg(G_h(t)) \equiv \alpha (\bmod~\pi)$
for a fixed value of $\alpha$.
The spacing of these points  is extremely 
regular, as determined by Stirling's formula.
Under the given hypotheses the perturbation coming
from the zeta function contribution is sufficiently
``small'' that the  
spacings of the perturbed zeros remain the same in the main
term of their asymptotics. The argument breaks down at
$h=0$, as it must, assuming the GUE Hypothesis is valid.

As a preliminary, we  collect needed estimates
on $\frac{\zeta'}{\zeta}(s)$. Let
$$
R_h(T) := \sup_{ T \le t \le T+1} 
|\frac{\zeta'}{\zeta}\left(\frac{1}{2} + h + it \right)|.
$$ 

\begin{lemma}~\label{le41} The function $R_h(T)$ satisfies
the following bounds.

(1) For $h > \frac{1}{2}$, there holds 
$$
R_h(T) = O (1),
$$
 where the $O$-constant depends on $h$.

(2) For $h = \frac{1}{2}$, and $|T| \ge 2$, 
$$
R_h(T) = O ( \frac{\log T}{\log\log T}).
$$ 

(3) Assuming the Riemann hypothesis, for $0 < h \le \frac{1}{2}$,
 $$
R_h(T) = O ( (\log T)^{1 - 2h}).
$$
where the $O$-constant depends on $h$.
\end{lemma}

\paragraph{Proof.}
(1) follows from logarithmic differentiation of the
Euler product for $\zeta(s)$, noting that $h > \frac{1}{2}$
is the region of absolute convergence.
(2) is shown in Titchmarsh \cite[5.17.4]{TH86}, in which 
the condition $|T| \ge 2$ is present to avoid the pole at $s=1$.
(3) is shown in Titchmarsh \cite[14.5.1]{TH86}.
~~~$\bsq$ \\

\paragraph{Proof of Theorem~\ref{th41}.}
Recall that 
$$
\varphi_h(t) = \arg( \xi(\frac{1}{2} + h + it)) = 
\Im \left(\log \xi(\frac{1}{2} + h + it) \right).
$$
and
$$
\tilde{\varphi}_h(t) =  \arg(G_h(t)) = 
\arg(\pi^{- (\frac{1}{4} + \frac{h}{2} + \frac{it}{2})} 
\Gamma(\frac{1}{4} + \frac{h}{2} + \frac{it}{2})).
$$
Under the hypotheses above (either (1) or (2)) 
both these functions are strictly
increasing functions of $t$. 
The proof of Theorem \ref{th31} gives
$$
\varphi_h(t)  = \tilde{\varphi}_h(t) + \pi + 
\arg(\zeta(\frac{1}{2} + h + it))
+ O (\frac{1}{|t| + 1}).
$$

From the formula \eqn{404}
for the argument of the Gamma factor
we deduce, for $0 \le \agh  \le 1,$ and $t \ge 0$ that
\begin{eqnarray}~\label{404a}
\tilde{\varphi}_h (t+ \agh) - \tilde{\varphi}_h (t) &=&
\frac{\agh}{2} \log \frac{t + \agh}{2} + \frac{t}{2} 
\log (1 + \frac{\agh}{2t})
+ \frac{\agh}{2} \log \frac{\pi}{e} + O(\frac{1}{|t|+1}) \nonumber \\
& = & \frac{\agh}{2} \log t + 
\frac{\agh}{2} \left( \frac{1}{2} - \log \frac{\pi}{2e}\right) +
O(\frac{\agh^2+1}{|t|+1}) \nonumber  \\
&=& \frac{\agh}{2} \log t + O ( \agh + \frac{1}{|t|+1}).
\end{eqnarray}
We now  define $\tilde{\gamma}_n(G_h)$ by
$$
\tilde{\varphi}_h \left(  \tilde{\gamma}_n(G_h) \right) =
\frac{\pi}{2} + n \pi.
$$
From the regular variation of this argument we can deduce for $t \ge 1$
that
\beql{405}
\tilde{\gamma}_{n+1}(G_h) - \tilde{\gamma}_{n}(G_h) =
\frac{2 \pi} {\log t} + O(\frac{1}{(\log t)^2}),
\eeq
although we will not directly use this in what follows.

The proof of Theorem~\ref{th31} deduces from \eqn{406} and \eqn{407} that
$$
\varphi_h(t)  = \tilde{\varphi}_h(t) + \pi + \arg(\zeta(\frac{1}{2} + h + it))
+ O (\frac{1}{|t| + 1}).
$$
We compare these quantities at $t$ and $t+\agh$, with $0 \le \agh \le 1$,
$t \ge 1$ to obtain
$$
\varphi_h(t + \agh) - \varphi_h(t) =
\left( \tilde{\varphi}_h(t+\agh) - \tilde{\varphi}_h(t) \right)
+\Delta \arg \left( \zeta(\frac{1}{2} + h + iu)\right)\Bigr|_{u=t}^{u= t + \agh} +
O \left( \frac{1}{|t| +1}\right).
$$
The change in argument of $\zeta(s)$ is estimated by
\begin{eqnarray*}
\Delta \arg \left( \zeta(\frac{1}{2} + h + iu) \right)|_{u=t}^{u= t + \agh}
& \le & | \int_{t}^{t+\agh} \Re \left( \frac{\zeta'}{\zeta}(\frac{1}{2} + h + iu)\right) du |\\
& \le & \int_{t}^{t+\agh} |\frac{\zeta'}{\zeta}\left(\frac{1}{2} + h + iu
\right)|du \\
& = & O(\agh R_h(t)).
\end{eqnarray*}
We bound the term $ R_h(t)$ using Lemma~\ref{le41}. 
Using \eqn{404a} we obtain for $0 < \agh \le 1$ that 
\beql{412}
\varphi_h(t + \agh) - \varphi_h(t) = \frac{\agh}{2} \log t + \tilde{R}_h(t) +
O(\frac{1}{|t| + 1}), 
\eeq
in which  the remainder term $\tilde{R}_h(t)$ satisfies
\begin{eqnarray*}
\tilde{R}_h(t) &= &~ O(\agh) ~~\mbox{if}~~ |h| > \frac{1}{2} \\
&= & ~ O( \frac{\agh \log t}{\log\log t}) ~~\mbox{if}~~ |h| = \frac{1}{2} \\
&= & ~ O( \agh (\log t)^{1- 2|h|}) 
~~\mbox{if RH holds and}~~0 < |h| \le \frac{1}{2}.               
\end{eqnarray*}
The global spacing in Theorem~\ref{th31} gives
$$\gamma_n(A_{h, \theta}) = 2 \pi \frac{n}{\log n} + 
O \left( \frac{n}{(\log n)^2}\right).
$$
We can now invert \eqn{412} to infer that, for 
$t = \gamma_n(A_{h, \theta})/ \log \gamma_n(A_{h, \theta})$, that 
$$
{\gamma}_{n+1}(A_{h, \theta}) - {\gamma}_{n}(A_{h, \theta}) =
\frac{2 \pi} {\log n} + O(\sR_h( \frac{n}{\log n}))
$$
in which the remainder term $\sR_h(t)$ satisfies
\begin{eqnarray*}
\sR_h(t) &= & ~ O(\frac{1}{ (\log t)^2} ) ~~\mbox{if}~~ |h| > \frac{1}{2} \\
&= & ~ O( \frac{1}{\log\log t}) ~~\mbox{if}~~ |h| = \frac{1}{2} \\
&= & ~ O( \frac{1}{(\log t)^{1+ 2|h|} }) 
~~\mbox{if RH holds and}~~0 < |h| \le \frac{1}{2}.               
\end{eqnarray*}
Now we know that
$$
\gamma_n(A_{h, \theta}) = \frac{1}{2\pi} \frac{n}{\log n} + 
O\left( \frac{n}{(\log n)^2} \right),
$$
so we conclude that the normalized zero spacings 
$$
\delta_n(A_{h, \theta}) := 
({\gamma}_{n+1}(A_{h,\theta}) - {\gamma}_{n}(A_{h, \theta})) \frac{1}{2\pi}
\log \frac{\gamma_n(A_{h, \theta})}{2\pi} 
$$
satisfy
$$
\delta_n(A_{h, \theta}) = 1 + 
O \left( (\log n) \sR_h( \frac{n}{\log n}) \right).
$$
The bounds on $\sR_h(t)$ give in all cases
$$ 
\delta_n(A_h)= 1 + O \left( \frac{1}{\log\log n} \right),
$$
under the stated hypotheses. This then gives the result for
$k$ consecutive spacings for any fixed value of $k \ge 1$. 
The same results 
applies to $\delta_n(B_{h, \theta})$ by an identical argument.
~~~$\bsq$ \\

%
%
%
%

\section{Differenced  $L$-functions}

The results of this paper extend without essential
change to all $GL(1)$ $L$-functions over $\QQ$. These
are the completed Dirichlet $L$-functions
$\xi(s, \chi)$, 
associated to  a primitive character
$\chi$ of conductor $N$, given by
$$
\xi(s, \chi) := (\frac{\pi}{N})^{- \frac{s+k}{2}}
\Gamma( \frac{s+k}{2}) L(\chi, s)
$$
where $k=0$ or $1$ according as $\chi(-1) = \pm 1,$
and $L(s, \chi)$ denotes the usual Dirichlet $L$-function.
These functions satisfy the 
functional equation 
$$
\xi(s, \chi) = \xi(1-s, \bar{\chi}),
$$
and are real-valued on the critical line. In this section
we formulate the analogous results, omitting detailed proofs.

\begin{lemma}~\label{le51}
(1) For $|h| \ge \frac{1}{2}$,
the function $E_h(s, \chi) := \xi(s+h, \chi)$ satisfies
\beql{503}
|E_h(s, \chi)| > |E_h( 1 - \bar{s}, \chi)|~~~~\mbox{for}~~
\Re(s) > \frac{1}{2}.
\eeq
 
(2) Assuming  the Riemann hypothesis holds for $\xi(s, \chi)$,
then the inequality  \eqn{503} is valid  for 
all nonzero $h$. 
\end{lemma}

\paragraph{Proof.}
The proofs are similar to those for Lemma~\ref{le21}.
We need to know
that $\xi(s, \chi)$    has a modified Hadamard product of the form 
\beql{504}
\xi(s, \chi) = e^{ A^{\ast}(\chi)+  B^{\ast}(\chi)s}
\prod_{\rho, \ast} \left( 1 - \frac{s}{\rho}\right),        
\eeq
where $\rho$ runs over the zeros of $L(\chi, s)$ inside the
open critical strip $0 < \Re(s) < 1$, and the
product is interpreted as the 
limit as $T \to \infty$
over all zeros with $|\rho| \le T$, with 
$$
B^{\ast}(\chi) := B(\chi) + \lim_{T \to \infty}
\sum_{|\rho| < T} \frac{1}{\rho}.
$$
This can be derived from the Hadamard product formula, using
the standard  asymptotic formula for zeros up to height $T$,
as given in Davenport~\cite[Chap. 16]{Da80}, to show that
the zeros inside  a box of side $1$ at heights $T$ and $-T$ respectively
combine to give a convergent sum.
In addition we have
\beql{505}
\Re (B^{\ast}(\chi)) = 0,
\eeq
which follows from  Davenport~\cite[p. 83]{Da80}, who shows that
$$ 
\Re (B(\chi)) = - \frac{1}{2} \sum_{\rho} (\frac{1}{\rho} +
 \frac{1}{\bar{\rho}}) = - \sum_{\rho} \Re(\frac{1}{\rho}).
$$
In particular  $\Re (B^{\ast}(\chi)) \ge 0$, and the rest of the
proof follows that of  Lemma~\ref{le21}.
~~~$\bsq$ \\

We now set
$$
E_{h, \theta}(s, \chi) := e^{i \theta} \xi(s+h, \chi) = 
A_{h, \theta}(s, \chi) -i B_{h, \theta}(s, \chi)
$$
with $ A_{h, \theta}(s, \chi)-i B_{h, \theta}(s, \chi)$ 
given as in Lemma~\ref{le22}.

\begin{theorem}~\label{th51}
Let $\chi$ be a primitive Dirichlet character that is 
non-principal, i.e. $\chi \ne \chi_0$.

(1) For $|h| \ge \frac{1}{2}$ and any $0 \le \theta < 2 \pi$,
the entire functions
$A_{h, \theta}(s, \chi)$ and $B_{h, \theta}(s, \chi)$ have 
 all their zeros on the critical line $\Re(s)= \frac{1}{2}.$
These zeros are all simple zeros, and they interlace.

(2) Assuming the Riemann hypothesis  for $L(s, \chi)$,
then  for $0 < |h| < \frac{1}{2}$ and any $0 \le \theta < 2 \pi$
the functions
$A_{h, \theta}(s, \chi)$ and $B_{h, \theta }(s, \chi)$ have 
have all their zeros on the critical line $\Re(s)= \frac{1}{2}.$
These zeros are all simple zeros, and they interlace.
\end{theorem}

\paragraph{Proof.}
These results are established by essentially the same proof as
Theorem~\ref{th21}, using Lemma~\ref{le51} in place of Lemma~\ref{le21}.
~~~$\bsq$ \\

Next we observe that an  analogue of Theorem~\ref{th31} holds, with
a similar proof, following Davenport~\cite[Chap. 16]{Da80}. 
It asserts that, for $|T| \ge 2$ ,
$$
N(T, A_{h, \theta}(s, \chi)) = \frac{1}{\pi} T \log T - 
\frac{1}{\pi}(\log (\frac{2\pi}{N}) + 1)T + 
O( \log T + \log N).
$$
It is valid unconditionally for $|h| \ge \frac{1}{2}$
and, under the Riemann hypothesis for $L(s, \chi)$, for
all nonzero $h$. Finally, using this asymptotic
formula,  one can establish by
a proof similar to Theorem~\ref{th41} the
following result.

\begin{theorem}~\label{th53}
Let $\chi$ be a primitive Dirichlet character which
is non-principal, and let $k \ge 1$ be given.

(1) For $|h| \ge \frac{1}{2}$ and $0 \le \theta < 2 \pi$
the functions
$A_{h, \theta}(s, \chi)$ and $B_{h, \theta}(s, \chi)$ 
both have limiting distributions of their $k$ consecutive
normalized zero spacings. This distribution
of normalized zero spacings consists  of a delta function
located at $x=1$, the ``trivial'' distribution.

(2) Assuming the Riemann hypothesis for $L(s, \chi)$, the same
results also hold for $0 < |h| < \frac{1}{2}$ and
$0 \le \theta < 2 \pi$, 
with the limiting distribution of $k$ consecutive normalized zero 
spacings existing and being
the ``trivial'' distribution.
\end{theorem}

\paragraph{Proof.} This is established in a manner similar to
Theorem~\ref{th41}. The key ingredient is an analogue of
the three parts of Lemma~\ref{le41}, which
cover $|h| > \frac{1}{2}$, $|h| = \frac{1}{2}$
and $0 < |h| < \frac{1}{2}$,
respectively.  Here the analogue of part (1),
that $R_h(T) = O(1)$ follows
from the absolute convergence of the Dirichlet series for
$-\frac{L'}{L}(s, \chi)$ in $\Re(s) >1$.  
The analogue of part  (3), comprising the estimate
$$
|\frac{L'}{L}(s, \chi) | = O ( (\log |T|)^{2 -2 \sigma})
$$
valid for  $\Re(s) = \sigma > \frac{1}{2}$,
with $|T| \to \infty$,
follows from a result 
 of Iwaniec and Kowalski ~\cite[Theorem 5.17]{IK04}, noting
that the Ramanjuan-Petersson conjecture assumed in
that theorem holds for
Dirichlet $L$-functions.
The  analogue of 
part (2), that
$$
|\frac{L'}{L}(1+iT, \chi) | + O\left( \frac{\log |T|}{\log\log |T|}\right)
$$
holds for $|T| \ge 3$, is deducible from the formula
\cite[Prop. 5.16]{IK04}, taking the test function
$\phi(y) = 1-y$ for $0 \le y \le 1$
and $\phi(y) =0$ for $y \ge 1$, and 
using the standard zero-free region 
excluding zeros when $\sigma > 1 - \frac{c}{\log |T|}$,
(\cite[Theorem 5.10]{IK04}), to bound  the 
contribution over zeros.
~~~$\bsq$
%
%
%
%
\section{Hilbert Spaces of Entire Functions}

We now consider the functions  
$A_{h, \theta}(s)$ and $B_{h, \theta}(s)$ from
the viewpoint of  the
de Branges theory of Hilbert spaces of entire functions
(\cite{deB68}).

We recall some fundamental facts about the 
de Branges theory of Hilbert spaces of entire
functions, following  \cite{La06},
which gives a translation into operator-theoretic
language of results stated in de Branges's book  \cite{deB68}.
 A {\em de Branges
structure function} $E(z)$ is any entire function
having the property that
\beql{701}
|E(z)| > |\overline{E(\bar{z})}| ~~~\mbox{for}~~~\Im(z) > 0.
\eeq
Associated to any
entire function $E(z)$, whether it satisfies
\eqn{701} or not, is a unique decomposition 
\beql{702}
E(z) = A(z) - i B(z)
\eeq
in which $A(z), B(z)$ are entire functions that are
real on the real axis. We define
$$
E^{\#}(z):= \overline{E(\bar{z})},
$$
 which is itself an entire function, then
$A(z) = \frac{1}{2}(E(z) + E^{\#}(z))$,
$B(z) = -\frac{1}{2i}(E(z) - E^{\#}(z))$.
The usefulness  of \eqn{701} is that it implies that
the functions
$A(z)$ and $B(z)$ have only real zeros,
and that these zeros interlace. This was shown
by de Branges \cite[Lemma 5]{deB59}, 
and is the content of Lemma~\ref{le22} above, taking $s= \frac{1}{2} -iz$.

We assign to a structure function $E(z)$ a Hilbert space $\sH(E(z))$ of
entire functions, by a formulation given below,
whose Hermitian scalar product for admissible functions
$f(z), g(z)$ takes the form
\beql{703}
\langle f(z), g(z) \rangle_{\sH(E)} := \int_{-\infty}^{\infty}
\frac{f(x) \overline{g(x)}}{|E(x)|^2} dx.
\eeq
The functions belonging to $\sH(E(z))$
are exactly those entire functions which have a finite norm 
$||f||_{E}^2:= \langle f(z), f(z) \rangle_{\sH(E)}$ and whose growth
with respect to $E(z)$ is controlled in the upper half-plane
$\CC^{+} := \{ z: ~\Im(z) > 0\}$, and in the
lower half-plane $ \CC^{-} := \{ z: ~\Im(z) < 0\}$, as follows.
 We require that the two functions
$\frac{f(z)}{E(z)}$ and $\frac{ \overline{f(\bar{z})}}{E(z)}$
each be of bounded type and nonpositive mean type in  $\CC^{+}$. 
A function $h(z)$ is of {\em bounded type} if it can be
written as a quotient of two bounded analytic functions in $\CC^{+}$
and it is of nonpositive mean type if it grows no faster
than $e^{\epsilon y}$ as $y \to \infty$ on the imaginary axis 
$\{iy~: y >0\}$, for each $\epsilon > 0$. 
De Branges gives several different characterizations
of the functions belonging to $\sH(E)$ (e.g.  \cite[Theorem 20]{deB68}) 
and shows in particular that this Hilbert space is never
trivial; it is at least one-dimensional. 

Also associated to the space $\sH(E)$ is a
(generally unbounded)  multiplication
operator $(M_z, \sD_z)$ defined on the domain of all
functions $f(z) \in \sH(E)$ such that $zf(z) \in \sH(E)$.
The closure of $\sD_z$ is either all of $\sH(E)$ or
is a subspace of codimension one in $\sH(E)$; we call
the former case the ``dense'' case. In the ``dense'' case
the operator $M_z$ is a closed symmetric operator,
with deficiency indices $(1,1)$ in von Neumann's sense.
In particular it possesses self-adjoint extensions,
the complete set of which form
a  one-parameter family identifiable
with $\bU(1)=\{e^{i \theta}: ~ 0 \le \theta < 2 \pi \}$. 
The structure function $E(z)$ can be viewed as singling out one such
self-adjoint extension, associated to the function $A(z)$,
as follows. As noted above, $A(z)$ only has real zeros.
This self-adjoint  extension has pure discrete simple
spectrum, given by the zeros of $A(z)$ (counted {\em without 
multiplicity}), with associated eigenfunction
\beql{704}
f_{\rho}(z) := \frac{A(z)}{z - \rho} \in \sH(E),
\eeq
and requiring that
\beql{705}
M_z( f_{\rho})(z) = \rho f_{\rho}(z)
\eeq
hold for all $\rho$. The domain of this
self-adjoint extension is $\sD_z(A) := 
\sD_z \oplus \CC[ f_{\rho}(z)]$ where we have adjoined
the eigenfunction for any  single zero; this domain
 is well-defined since
$$
f_{\rho}(z)- f_{\rho'}(z) = (\rho' - \rho) \frac{A(z)}{(z-\rho)(z- \rho')}
\in \sD_z.
$$
For self-adjoint operators with discrete spectrum
a standard fact is that the eigenfunctions $f_{\rho}(z)$ for different
zeros $\rho$ are orthogonal in the Hilbert space  
$ \sH(E)$.
The complete family of self-adjoint extensions of the operator
$M_z$ is given by the spaces $\sH(E_{\theta})$ where
$E_{\theta}(z) := e^{i \theta} E(z)$ for $0 \le \theta < 2 \pi$.
The Hilbert spaces $\sH(E_{\theta})$ are identical
(same functions, same scalar product), but the variation in
structure functions picks out different self-adjoint extensions,
using  the decomposition $E_{\theta}(z) = A_{\theta}(z) - i B_{\theta}(z)$.
Note that $B(z) = A_{\theta}(z)$ with $\theta= \frac{\pi}{2}$.
The de Branges theory thus gives a spectral interpretation of
the zeros of $A(z)$ (counted without multiplicity) as the
spectrum of the self-adjoint operator $(M_z, \sD_z(A))$ in $\sH(E(z)).$

The de Branges theory can be regarded as giving
a normal form for a particular kind
of  non-selfadjoint operator $(M_z, \sD_z)$ with  deficiency
indices $(1,1)$, much as the spectral theorem for 
(unbounded) self-adjoint operators gives
a normal form as a multiplication operator, acting on a Hilbert
space given with a spectral measure, and with a specified domain.
The structure function $E(z)$ then encodes the data analogous to
the spectral measure.
The major content of de Branges theory is  a Fourier-like transform
(which might be called the {\em de Branges transform})
that converts this multiplication operator into a $2 \times 2$
matrix  integral
operator of a certain kind acting on a different
Hilbert space, much as the Fourier transform
takes the multiplication operator to a differential operator.
We do not address that aspect of the de Branges theory here.

The results of \S2 can be reformulated in 
 terms of the de Branges theory as follows.

\begin{lemma}~\label{le71}
(i) For $h \ge \frac{1}{2}$ the function
\beql{710}
E_h(z) := \xi( \frac{1}{2} + h - iz)
\eeq
is a de Branges  structure function, i.e. $|E_h(z)| > |E_h(\bar{z})|$
when  $\Im(z) > 0$.

(ii) If the Riemann hypothesis holds, then for all $h \ne 0$
the function $E_h(z)$ is a de Branges structure function.
\end{lemma}

\paragraph{Proof.} This is a restatement of Lemma~\ref{le21}.~~~$\bsq$ \\

Under the hypotheses of Lemma~\ref{le71} it follows 
that $E_{h, \theta}(z):=e^{i \theta}E_h(z)$ 
is a de Branges structure
function, whose associated decomposition \eqn{702} is
$$
E_{h, \theta}(z) = \tilde{A}_{h,\theta}(z) - i \tilde{B}_{h, \theta}(z)
$$
in which  
$\tilde{A}_{h,\theta}(z) := A_{h, \theta}(\frac{1}{2} -iz)$ and   
$\tilde{B}_{h,\theta}(z) := B_{h, \theta}(\frac{1}{2} -iz)$.
The de Branges theory now applies to give
a ``spectral'' interpretation of the zeros of $A_{h, \theta}(s)$,
or, equivalently, of the zeros of 
$$
f_{h, \theta}(t) := \Re \left( e^{i \theta} \xi(\frac{1}{2} + h + it)\right).
$$
as the eigenvalues of the self-adjoint operator 
$(M_{z}, \sD_z(A_{h, \theta}))$ acting in the de Branges Hilbert space
$\sH(E_{h, \theta}(z))$.
Theorem~\ref{th21} established that 
the zeros of $A_{h, \theta}(s)$  are  simple zeros for $|h| \ge \frac{1}{2}$,
and, assuming the Riemann hypothesis, for all nonzero $h$.

The particular de Branges space with structure function
$E(z) = \xi(1-iz)$ was considered  by de Branges \cite{deB86}
as a possible approach to the Riemann hypothesis.
This structure function is included in the family 
$$
E_{h, \theta}(z)= A_{h,\theta}(\frac{1}{2}-iz) - 
i B_{h,\theta}(\frac{1}{2} - iz)
$$ 
of Theorem~\ref{th21}, on taking $h = \frac{1}{2}$, $\theta=0$.
This result supplies a
 proof that $\xi(1-iz)$ is a structure function.
In ~\cite{deB86} and  \cite{deB92}
de Branges proved  general theorems giving  sufficient
conditions on the inner product of a Hilbert space of 
 entire functions $\sH(E(z))$ which imply that the associated
(normalized) structure function $E(z)$ necessarily has  all its zeros on
the line $\Im(z) = - \frac{1}{2}$.
If any of these theorems applied to the 
de Branges space $\sH(\xi(1 - iz))$,
the Riemann hypothesis would follow. 
However  Conrey and Li~\cite{CL00}
recently showed that
the hypotheses of these theorems are
not satisfied for 
the de Branges Hilbert spaces $\sH(E(z))$
with $E(z) = \xi_{\chi}(1 - iz)$ with  $\chi$
either the trivial character $\chi_0$ or
 the real character $\chi_{-4}$.
There do exist de Branges spaces satisfying
the inner product conditions of \cite{deB86}, 
presented in  Li \cite{Li00}, see also \cite{Li97}.

In \cite{La06} we
formulate another connection between the Riemann
hypothesis for Dirichlet $L$-functions and certain 
Hilbert spaces of entire
functions. This connection is conditional,
in that the  associated de Branges space exists
if and only if the Riemann hypothesis holds for the
corresponding $L$-function.

%
%
%
%
\section{Concluding Remarks}

(1) This paper studied
a two-parameter deformation $A_{h, \theta}(s)$
of the Riemann $\xi$-function
using the  parameters $(h, \theta)$. It showed that
under the RH this deformation
 preserves a ``Riemann hypothesis'' condition,
and proves this holds unconditionally when  
$|h| \ge \frac{1}{2},$ for all $\theta$. 
The deformations  $(h, \theta)$ 
preserve a ``functional equation'',
namely
$$
A_{h, \theta}(s)= \overline{A_{h, \theta}(1- \bar{s}}),
$$
which encodes the property that  $A_{h, \theta}(s)$
is real on the critical line,
embodying the reflection principle. The deformations  also
preserve a second `` functional equation''
$$
B_{h, \theta}(s)= \overline{B_{h, \theta}(1- \bar{s}}).
$$
However we do not know of any  analogue of an
Euler product (or Hecke operator factorization)
that is preserved under such deformations.

(2) The arguments of \S5 
extend with little change to appropriate 
automorphic $L$-functions over  $GL(N)$
(principal $L$-functions over $GL(N)$ over the
rational field $\QQ$) for 
all $N \ge 2$, under suitable extra  hypotheses.
indicated below.
Given an $L$-function, one  constructs an analogue
 of the $\xi$ function,
adjusted to be real on the critical line,
as is done in 
\cite[\S2]{La04b}. 
One has global  zero density estimates that the
number of zeros to height $T$ grows like 
$\frac{N}{\pi}T \log T$, see  
Iwaniec and  Kowalski \cite[Theorem 5.8]{IK04}, 
or Lagarias \cite[Theorem 2.1]{La04b}.
One can prove an unconditional result for 
$|h| > \frac{1}{2}$, since the 
analogue of Lemma \ref{le41}(1) holds: the
 Dirichlet series of $L(s, \pi)$ converges
absolutely for $\Re(s) > 1$, and this implies
an $O(1)$ bound for the analogue
of $R_h(T, \pi)$ where the $O$-constant depends on $h$.
The convergence of the
Dirichlet series is formulated in  Lagarias \cite[Theorem 2.1]{La04b},
a result  that goes back to work of Jacquet and Shalika.
A conditional result can also
be proved, valid for all nonzero $h$,
assuming the truth of both  the 
Riemann hypothesis and the the Ramanujan-Petersson conjecture
for the given 
given automorphic $L$-function attached to an irreducible, cuspidal
unitary representation of $GL(N)$ over the adeles. 
(The Ramanujan-Petersson conjecture states that the 
local parameters  $|\alpha_i(p)|$ are of absolute value 
at most one, see \cite[p. 95]{IK04}).
Under
these hypotheses a suitable analogue of Lemma \ref{le41}(3) holds 
using Theorem 5.17  of Iwaniec and Kowalski ~\cite{IK04}.

(3)  Theorem~\ref{th41} implies  that  non-trivial zero
spacing statistics (GUE-like statistics) for these functions
are associated precisely  with
the critical line, i.e. $h=0$.  
It has been asserted  that the Riemann hypothesis, if true,
is ``just barely true''. One precise formulation  of
this has to do with the de Bruijn-Newman constant,
see Csordas, Smith and Varga \cite{CSV94}. 
Thus since $\xi(s)$ presumably
satisfies GUE, we would not expect $\xi(s)$ itself to be produced by
averaging shifts of another function with all its zeros on
the critical line. We therefore 
ask: Is it true that any entire function
$G(s)$ such that for some $h >0$, 
$$
\xi(s) = \frac{1}{2}\left( G(s + h) + G(s- h) \right)
$$
necessarily has the property that not all its zeros lie
on the critical line $\Re(s) = \frac{1}{2}$?

(4) The result of Theorem~\ref{th41} also illustrates
a heuristic principle that ``averaging'' and  ``differencing'' operations
smooth the spacings of zeros. 
The operator 
$$
T_{h}(f)(s) := \frac{1}{2}(f(s+h) + f(s-h))
$$
is  a convolution-type operator which convolves with 
 a discrete probability
measure (with masses of weight $\frac{1}{2}$ at
the points $\pm h$).
Such convolution operators 
may  be expected to smooth zeros under some conditions.
on the operator,
as studied in Cardon and Nielsen \cite{CN02}, for functions
in the Laguerre-Polya class.
In the case $f(s)= \xi(s)$ a heuristic goes as follows.
The zeros of a shifted function $\xi_h(s)$ ``feel the effect'' of
zeros of the unshifted function within distance $O(1)$.
Assuming the Riemann hypothesis, there
 is a regular asymptotics for zeros over
a range $[T, T+1]$ of this length for the $\xi$-function,
namely 
$$ N(T+1) - N(T) = \frac{1}{2\pi} \log T  + O(\frac{\log T}{\log\log T}),$$
see Titchmarsh~\cite[Theorem 14.13]{TH86}.
The heuristic is that the effect of averaging is to make the zeros
repel each other. There are definite restrictions needed to formulate
in a general theorem,
for example bounds on the growth rate (entire function of order less than
$2$) amd also on the location of zeros of the function $f(s)$.
D. Cardon (private communication)
has  examples showing that show this heuristic cannot be valid without
further hypotheses, e.g. that 
the density of zeros to height $T$ grow faster than linear in $T$.

(5) One might also study the distribution of normalized
zeros of the derivative $\xi'(s)$ of the $\xi$-function. Assuming the
Riemann hypothesis, all the zeros of $\xi'(s)$ lie on the 
critical line. If the zeros of $\xi(s)$ satisfy the
GUE hypothesis, it seems reasonable to expect that the 
normalized spacings of zeros of $\xi'(s)$
also have a limiting distribution.
This distribution cannot be the ``trivial'' distribution
concentrated at equal normalized spacings of size
$1$, because the zeros of $\xi'(s)$ interlace
with those of $\xi(s)$, and the GUE distribution predicts a
positive probability of two consecutive normalized spacings each of which 
is at most $\frac{1}{2} - \delta$, for fixed positive $\delta$,  
so the associated normalized zero spacing of $\xi'(s)$ is 
with positive probability at most 
$1 - 2 \delta$. However this distribution is certainly not 
the GUE distribution. The distribution of
normalized zero spacings of  $\xi'(s)$ should
be a new distribution, whose form is expected to  be predictable using
random matrix theory analogues. It should be more concentrated
near the  unit spacing than the GUE distribution, by analogy with 
results of Farmer and Rhodes \cite{FR04}.

(6) It seems plausible  that there will be  
a random matrix theory analogue of Theorem~\ref{th41}.
This would concern the behavior of the roots of
differenced characteristic polynomials
of random unitary matrices drawn from $U(N)$, followed by
letting $N \to \infty.$ 
One may also expect there  to
be a random matrix theory analogue of the distribution of
the zeros of $\xi'(s)$. We hope to address these 
questions elsewhere.


\noindent{\it AMS 2000 Subject Classification:}  
11M26 (Primary); 11M41, 46C05, 47B32  (Secondary) \\

\noindent email: {\tt lagarias@umich.edu}
\end{document}